\documentstyle[12pt]{article}
\newcommand{\const}{\mathop{\rm const}\limits}

\newcommand{\Var}{\mathop{\rm Var}\limits}

\newcommand{\Law}{\mathop{\rm Law}\limits}

\newcommand{\Vol}{\mathop{\rm Vol}\limits}

\begin{document}

\begin{center}

{\bf  MONTE CARLO COMPUTATION OF MULTIPLE } \\

\vspace{3mm}

{\bf WEAK SINGULAR INTEGRALS}  \\

\vspace{3mm}

{\bf OF SPHERICAL AND VOLTERRA'S TYPE.} \par

\vspace{4mm}

 $ {\bf E.Ostrovsky^a, \ \ L.Sirota^b } $ \\

\vspace{4mm}

$ ^a $ Corresponding Author. Department of Mathematics and computer science, Bar-Ilan University, 84105, Ramat Gan, Israel.\\

\vspace{3mm}

E-mail: \ eugostrovsky@list.ru\\

\vspace{3mm}

$ ^b $  Department of Mathematics and computer science. Bar-Ilan University,
84105, Ramat Gan, Israel.\\

\vspace{3mm}

E-mail: \ sirota3@bezeqint.net \\

\vspace{4mm}
                    {\bf Abstract.}\\

 \end{center}

 \vspace{4mm}

 We offer a simple  method  Monte  Carlo  for computation of Volterra's and spherical type multiple integrals
with weak (integrable) singularities. \par

 An elimination of infinity of variance is achieved by incorporating singularities in the density, and we offer
 a highly effective way for generation of appeared multidimensional distribution.\par

 We extend offered method onto multiple Volterra's and spherical integrals with weak singularities containing parameter. \par

\vspace{4mm}

 \par

\vspace{4mm}

  \par

\vspace{4mm}

\section{Introduction. Notations. Statement of problem.}

\vspace{3mm}

 Let us denote by $  S_n  $ the $ n \ - $ dimensional polyhedron (simplex) of a form

$$
S(n) = \{ \vec{s} = (s_1, s_2, \ldots,s_n ): \ 0 < s_1 < s_2 < \ldots < s_n < 1 \}.
$$

 It is easy to calculate that $ \Vol(S(n)) = 1/n!. $\par
Let also $ \alpha = \vec{\alpha} = ( \alpha_1, \alpha_2, \ldots, \alpha_n )  $ be a constant vector  such that
$ 0 \le \alpha_k < 1, \ k = 1,2,\ldots,n. $ \par

{\bf We investigate in this article the problem of numerical computation by the Monte  Carlo method of the following
integral of Volterra's type} \hspace{5mm}

$$
I_{\alpha,n} [z] :=
$$

$$
 \int_0^{1} ds_n \int_0^{s_n} \ ds_{n-1} \int_0^{s_{n-1}} ds_{n-2} \ldots \int_0^{s_2} ds_1
\frac{z(s_1,s_2, \ldots, s_n) \ ds_1 ds_2 \ldots ds_n}{s_1^{\alpha_1} (s_2 - s_1)^{\alpha_2} (s_3 - s_2)^{\alpha_3} \ldots (s_n - s_{n-1})^{\alpha_n} } =
$$

$$
 \int \int \ldots \int_{S(n)}
\frac{z(s_1,s_2, \ldots, s_n) \ ds_1 ds_2 \ldots ds_n}{s_1^{\alpha_1} (s_2 - s_1)^{\alpha_2} (s_3 - s_2)^{\alpha_3} \ldots (s_n - s_{n-1})^{\alpha_n} },
 \ n = 2,3,\ldots; \eqno(1.1)
$$

$$
I_{\alpha,1}[z] := \int_0^1 \frac{z(s_1)}{s_1^{\alpha_1}} \ ds_1. \eqno(1.1a)
$$

\vspace{3mm}

This kind of such integrals appear in the reliability theory \cite{Grigorjeva2}; in the theory of integral equations of Volterra's type with
weak singularity in kernel \cite{Grigorjeva1}, \cite{Heinrich1}, \cite{Heinrich2}; in the investigation of the spectres of some integral
operators with singularities \cite{Ostrovsky1}, p. 310;  in the investigation of local times of random processes \cite{Molchanov1};
in the numerical analysis  of the Navier - Stokes equation
\cite{Giga1}, \cite{Giga4}, \cite{Kato1}, \cite{Kato2}, \cite{Ostrovsky111}, \cite{Ostrovsky201} etc. \par

 Let us discuss the last thesis  in more detail.  The so-called mild solution $ u = u(x,t) $ of the Navier-Stokes equation in the whole space
$  x \in R^d  $ during  its lifetime $ t \in [0,T], \ 0 < T = \const \le \infty $ may be represented under some simple conditions as a limit as
$ m \to \infty, m = 0,1,2, \ldots $ of the following recursion:

 $$
 u_{m+1}(x,t) = u_0(x,t) + G[u_m, u_m](x,t), m=1,2,\ldots
 $$
where $ u_0(x,t) $ is the solution of heat equation  with correspondent initial value and right-hand side
and $ G[u,v] $ is bilinear unbounded pseudo-differential operator, \cite{Kato1}.\par
 The detail expression for each iteration starting at the value $ m = 2 $ contains multiple integrals of a form (1.1) with the values
 $ \alpha_k = 1/2, \ k = 1,2,\ldots,n   $ and hence  may be computed  for example by means of the Monte - Carlo method,  in addition to
 the deterministic methods, see \cite{Temam1}, \cite{Vishik1}.  \par
Note that the second iteration is considered in the article \cite{German1}. \par

\vspace{3mm}

 The case of spherical multiple integrals with weak singularities alike (1.1), i.e. when instead the polygon $  S(n) $
states (unit) ball    $ B = B_1 $ will be considered further. \par

 \vspace{3mm}

\section{ The essence of the method.} \par

\vspace{3mm}

 Let us denote for brevity for the values $ s = \vec{s} \in S(n) $

$$
R_{\alpha,n,S}(s) = R_{\alpha,S}(s)  = s_1^{-\alpha_1} (s_2 - s_1)^{-\alpha_2} (s_3 - s_2)^{-\alpha_3} \ldots (s_n - s_{n-1})^{-\alpha_n}; \eqno(2.1)
$$
then

$$
 I_{\alpha,n} [z] =  \int_{S(n)} z(s) \  R_{\alpha,S}(s) \ ds. \eqno(2.2)
$$

 The "direct" probabilistic representation for $ I_{\alpha,n}[z],  $  i.e. the expression of a form

 $$
 I_{\alpha,n}(z) = n! \cdot {\bf E} z(\eta) R_{\alpha,S}(\eta)
 $$
where the random vector $ \eta = \{\eta_1, \eta_2, \ldots, \eta_n \} $ has an uniform distribution in the simplex $  S(n) $
 leads in general case, more exactly, when $ \exists k: \alpha_k \ge 1/2, $
 to the possibility of infinite variance in the integrand  expression.
  And then, for the error estimate instead of the CLT  will have to apply the Stable Limit Theorem (SLT),
which drastically reduces the rate of convergence of the proposed method  \cite{Ostrovsky1}, chapter 5, section 5.14.  \par

 Recall that this is the case of the problem of the Navier-Stokes equation, in which all the values $ \alpha_k $
are equal exactly to 1/2.\par

\vspace{3mm}

  We offer  hence {\it another } probabilistic representation for $ I_{\alpha,n}[z]. $ \par

 Note first of all that

$$
\int \int \ldots \int_{S(n)}\frac{ds_1 ds_2 \ldots ds_n}{s_1^{\alpha_1} (s_2 - s_1)^{\alpha_2} (s_3 - s_2)^{\alpha_3} \ldots (s_n - s_{n-1})^{\alpha_n} } =
K_{n,S}(\alpha) = K_{n,S}(\vec{\alpha}), \eqno(2.3)
$$
where

$$
K_{n,S}(\alpha) = K_{n,S}(\vec{\alpha}) \stackrel{def}{=} \frac{\prod_{k=1}^n \Gamma(1-\alpha_k)}{\Gamma(1 + \sum_{k=1}^n (1 - \alpha_k))},\ 0 \le \alpha_k < 1,
\eqno(2.4)
$$
and $ \Gamma(\cdot) $ is ordinary Gamma function. \par

For example, let $  \alpha $ be arbitrary  number such that $ 0 \le \alpha < 1. $ Denote $ \beta = \beta(\alpha) =1 - \alpha. $
Then

$$
W_n(\beta) := K_{n,S}(\alpha,\alpha, \ldots, \alpha) = \int_{S(n)} R_{\alpha, \alpha, \ldots,\alpha; n,S}(s) \ ds =
$$

$$
\int \int \ldots \int_{S(n)}\frac{ds_1 ds_2 \ldots ds_n}{s_1^{\alpha} (s_2 - s_1)^{\alpha} (s_3 - s_2)^{\alpha} \ldots (s_n - s_{n-1})^{\alpha} } =
\frac{\Gamma^n(\beta)}{\Gamma(1 + n \beta)}.
$$

\vspace{3mm}

 Evidently, $ \lim_{n \to \infty} W_n(\beta) = 0.  $\par

 Note in addition  the continuity and monotonicity of the function  $ \beta \to W_n(\beta):  $

$$
0 < \beta_1 < \beta_2 \le 1 \ \Rightarrow W_n(\beta_1) > W_n(\beta_2).
$$
 In particular,

$$
W_n(1/2) = \frac{\pi^{n/2}}{ \Gamma(1 + n/2)},
$$

$$
\frac{1}{n!} = \frac{\Gamma^n(1)}{\Gamma(1 + n)} = W_n(1) < \frac{\Gamma^n(\beta)}{\Gamma(1 + n \beta)} = W_n(\beta), \ 0 < \beta < 1.
$$

\vspace{3mm}

{\it  The following function $  h_{\alpha}(s) =  h_{\vec{\alpha}}(s), \ s \in S(n),  $
  could be chosen as a density of distribution with support on the   simplex} $  S(n): $

$$
h_{\alpha}(s) = \frac{R_{\alpha,n}(s)}{K_{n,S}(\alpha)}.\eqno(2.5)
$$

\vspace{3mm}

{\bf Definition 2.1.}  The random vector $  \kappa = \kappa_{\alpha,n} = \vec{\kappa}  = \vec{\kappa}_{\alpha,n}   $ with values in
the polygon $  S(n) $ has by definition a {\it polygonal Beta distribution,} write:   $  \Law(\kappa) = PB(\alpha,n), $   iff it has a
density $  h_{\alpha}(s), \ s \in S(n).  $\par

 On the other word,

$$
{\bf P}(\kappa \in G) = \int_G h_{\alpha}(s) \ ds \stackrel{def}{=} \mu_{\alpha,n}(G), \ G \subset S(n). \eqno(2.6)
$$

 Evidently, $ \mu_{\alpha,n}(\cdot) $ is the probabilistic Borelian measure on the set $  S(n).  $ \\

\vspace{3mm}

 The expression for the source integral (1.1) may be represented as follows.

$$
I_{\alpha,n}[z] = K_{n,S}(\alpha) \cdot \int_{S(n)} \frac{z(s) \ R_{\alpha,n}(s) \ ds}{ K_{n,S}(\alpha)} =
$$

$$
 K_{n,S}(\alpha) \cdot \int_{S(n)} z(s) h_{\alpha}(s) \ ds = K_{n,S}(\alpha) \cdot {\bf E} z(\kappa),\eqno(2.7)
$$
where the random vector $   \kappa $ has the polygonal Beta distribution:   $  \Law(\kappa) = PB(\alpha,n). $ \par

\vspace{3mm}

 Let us estimate the second moment of the variable $  K_{n,S}(\alpha) \cdot z(\kappa):  $

$$
{\bf E} \left[ K_{n,S}(\alpha) \cdot z(\kappa) \right]^2 = K^2_{n,S}(\alpha) \cdot  ||z||^2L_2(\mu_{\alpha,n}). \eqno(2.8)
$$

 The last expression is finite if for instance the function $ z(\cdot) $ is bounded. \par
  Moreover, it is interest to note that the variance of the r.v.  $  K_{n,S}(\alpha) \cdot z(\kappa),  $ i.e. the expression

$$
\sigma^2_{n,S}[z] := \Var \left\{ K_{n,S}(\alpha) \cdot z(\kappa) \right\} = K^2_{n,S}(\alpha) \cdot  ||z||^2L_2(\mu_{\alpha,n}) -
I^2_{\alpha,n}[z]
$$
tends very  rapidly to zero as $  n \to \infty,  $  if for example $  \sup_k \alpha_k < 1.  $ \par
  Thus, we can use by computation of this integral and by error estimation  the classical Monte -  Carlo method with application of
Central Limit Theorem (CLT):

$$
I^{(N)}_{\alpha,n}[z] := N^{-1} \sum_{i=1}^N K_{n,S}(\alpha) \ z(\kappa_i), \eqno(2.9)
$$
where $ \kappa_i $ are independent random vectors  with polygonal beta distribution.  The variance $ \sigma^2_{n,S}[z]  $ may be
estimated as well as the integral $ I_{\alpha,n}[z]. $\par
 Note that the "regular" case, i.e. when $ \vec{\alpha} = 0 $ may be obtained by substitution $ \beta = 1. $ For example,

$$
W_n(\beta(0)) = W_n(1) = 1/n!,
$$

$$
{\bf E} \left[ K_{n,S}(0) \ z(\kappa) \right]^2 = (n!)^{-2} \cdot ||z||^2L_2(\mu_{0,n}).
$$

Mentioned here  method is  called in the theory of Monte Carlo "incorporating singularities in the density."\\

\vspace{3mm}

\section{ Generation of used random vectors.   }

\vspace{3mm}

  Let us dwell briefly on the issue of the generation of this distribution on the basis of the standard uniform
$   U(0,1) $ generator. Our purpose in this section is to prove that the generation of multidimensional polygonal beta
distribution is not much harder as in one dimensional case.\par

 \vspace{3mm}

 {\bf 1.} \ The   one-dimensional case $  n = 1 $ is very simple, as long as the r.v. $  \kappa_1 $ has a power distribution
 with the density:

$$
f_{\kappa_1}(x) = (1 - \alpha_1) x^{-\alpha_1}, \ x \in (0,1).
$$

 The (cumulative) distribution function $ F_{\kappa_1}(x)  $ has a view:

$$
F_{\kappa_1}(x) = (1 - \alpha_1) \int_0^x z^{-\alpha_1} \ dz  = x^{1 - \alpha_1}, \ x \in [0,1];
$$
hence

$$
F_{\kappa_1}^{-1} (x) = x^{1/( 1 - \alpha_1 ) }, \ x \in [0,1].
$$
Thus, the power distribution may be generated very simple by means of the inverse function method.\\

\vspace{3mm}

{\bf 2. } A two-dimensional case $  n = 2. $ \par

 Let the two-dimensional random vector $ (\xi, \eta) $ be polygonal beta distributed:

$$
f_{\xi,\eta} (x,y) = \frac{1}{K_{2, S(2)}(\alpha_1,\alpha_2)} \ \frac{I(0 < x < y < 1)}{x^{\alpha_1} \ (y - x)^{ \alpha_2 } },
$$
 where

 $$
 0 \le \alpha_1, \alpha_2 < 1; \hspace{6mm} I(A) = 1, (x,y) \in A, \ I(A) = 0, (x,y) \notin A.
 $$

We calculate the {\it marginal} density  $ f_{\eta}(y), \ y \in (0,1): $

$$
f_{\eta}(y) = \frac{1}{K_{2, S(2)}(\alpha_1,\alpha_2)} \ \int_0^y \frac{dx}{x^{\alpha_1} \ (y-x)^{\alpha_2}  } =
(2 - \alpha_1 - \alpha_2) \  y^{1 - \alpha_1 - \alpha_2  },
$$
i.e. the r.v. $ \eta $ has the power distribution with parameter $ 1 - \alpha_1 - \alpha_2. $\par
As for the {\it conditional} density of distribution $ f_{\xi}(x/\eta = y), $ we have:

$$
f_{\xi}(x/\eta = y) = C(\alpha_1, \alpha_2) \ \frac{I(0 < x < y < 1)}{x^{\alpha_1} \ (y - x)^{\alpha_2}}.
$$
 The last relation implies that the r.v. $ \xi $ has under condition $ \eta = y $ the well-known one-dimensional
beta distribution on the interval $ (0,y) $  with parameters $  (\alpha_1, \alpha_2). $

\vspace{3mm}

 Note that in the computer system MATLAB there is a command $  R = betarnd(A,B) $ which generated a (sequence) of the independent one-dimensional
 beta distributed random (pseudo-random) variables on the set  $  (0,1) $ with parameters $ A,B. $ \par

 \vspace{3mm}

{\bf n.} Multidimensional case. \par

 Let the random vector $ \kappa = \vec{\kappa} = (\kappa_1, \kappa_2, \ldots, \kappa_n) $ has a polygonal beta distribution.
It is easy  to calculate

$$
f_{\kappa_2, \kappa_3, \ldots, \kappa_n}(x_2, x_3, \ldots, x_n)  =
\int f_{ \kappa_1,\kappa_2, \kappa_3, \ldots, \kappa_n}(x_1,x_2, x_3, \ldots, x_n) dx_1 =
$$

$$
 C \cdot \frac{I(0 < x_2 < x_3 \ldots < x_n < 1)}{x_2^{1 - \alpha_1 - \alpha_2} (x_3 - x_2)^{\alpha_3} \ldots (x_n - x_{n - 1})^{\alpha_n} },
$$
i.e. the $  (n-1) $ dimensional subvector $ (\kappa_2, \kappa_3, \ldots, \kappa_n ) $  has also polygonal beta distribution.\par
 Thus, the problem of $  n   $ dimensional polygonal beta distribution random generating may be easy reduced to the $ (n-1)  $ dimensional. \par

\vspace{3mm}

 {\bf Example.}  Let us consider  an important two-dimensional case $  n = 2 $ with $ \alpha_1 = \alpha_2 = 1/2; $
 then the correspondent  one-dimensional function of distribution of a second component $  \kappa_2 $ has a form

$$
F_{\kappa_2} (z) = \pi^{-1} \int_0^z \frac{dx}{\sqrt{x(1-x)}} = \pi^{-1} [\arcsin(2z - 1) + \pi/2   ], \ z \in [0,1].
$$

 The inversion function has an explicit view:

 $$
 F_{\kappa_2}^{-1}(z) = 0.5 + 0.5 \sin (\pi(z-1/2)), \ z \in [0,1].
 $$

\vspace{4mm}

\section{Spherical case}

\vspace{4mm}

 Let us consider the following  example. Let $ A = (A_1,A_2, \ldots, A_n) $ be a $ n \ -   $ dimensional numerical vector such that
$   -1 < A_i \le 0  $  and $ D := \sum_{k=1}^n A_k + n > 0. $\par
For all the $ n \ -  $ dimensional numerical vector $ x = (x_1, x_2, \ldots, x_n) \in R^n $ we define as
usually  the {\it monomial}

$$
x^A = |x_1|^{A_1} |x_2|^{A_2} \ldots |x_n|^{A_n}.
$$
 The unit ball with the center in origin in the classical Euclidean distance will be denoted by $  B = B_1: $

$$
B = B_1 = \{x, \ x \in R^n, \ |x| \stackrel{def}{=} \sqrt{(x,x)}\le 1 \}.
$$

 It is known, see e.g.  \cite{Cabre1}

$$
\int_B x^A dx =  K_{n,B}(\alpha) = K_{n,B}(\vec{\alpha} ) \stackrel{def}{=} \frac{\prod_{k=1}^n  \Gamma( (A_k + 1)/2  )}{ \Gamma(D/2 + 1)}.
\eqno(4.1)
$$

 For example, if $ \alpha_1 = \alpha_2 = \ldots = \alpha_n = 2 \gamma -1, \ \gamma = \const > 0,   $ then

$$
K_{n,B}(\alpha, \alpha, \ldots, \alpha)  = \frac{\Gamma^n(\gamma)}{\Gamma(1 + n \gamma)} = W_n(\gamma).
$$
The last expression tends again to zero  as $  n \to \infty $  very rapidly, as well as $ W_n(\beta). $\par

\vspace{3mm}

{\bf Definition 4.1.} The $ n \ - $ dimensional random vector $  \zeta  = \zeta_{\alpha,n} $ with values in the unit ball
 $  B  $ has a {\it ball beta distribution} $ BB(\vec{\alpha}) = BB(\alpha), $ if it has a density  $  g_{\alpha}(x), \ x \in B $ of a form

$$
g_{\alpha}(x) \stackrel{def}{=} \frac{1}{K_{n,B}(\alpha) }\cdot  x^A: \hspace{4mm} {\bf P}(\zeta \in G) = \int_G g_{\alpha}(x) \ dx
\stackrel{def}{=} \nu_{\alpha,n}(G), \ G \subset B. \eqno(4.2)
$$

Obviously, $ \nu_{\alpha,n}(\cdot) $ is probabilistic Borelian measure on the set $  B.  $ \\

\vspace{3mm}

 {\it We consider in this section the problem of Monte Carlo computation of the multiple integral}

$$
J_{\alpha,n}[z] \stackrel{def}{=} \int_B |x|^A \ z(x) \ dx. \eqno(4.3)
$$

 As before, the "direct"  simulation leads in general case to the infinity of variance, therefore we need to transform this integral:

$$
J_{\alpha,n}[z] = K_{n,B}(\alpha) \cdot \int_B z(x) \ g_{\alpha}(x) \ dx = {\bf E} [ K_{n,B}(\alpha) \cdot z(\zeta)], \eqno(4.4)
$$
where the random vector $ \zeta $ has the ball beta distribution  $ BB(\alpha).  $\par

\vspace{3mm}

 Let us estimate the second moment of the variable $  K_{n,B}(\alpha) \cdot z(\zeta):  $

$$
{\bf E} \left[ K_{n,B}(\alpha) \cdot z(\zeta) \right]^2 = K^2_{n,B}(\alpha) \cdot  ||z||^2L_2(\nu_{\alpha,n}). \eqno(4.5)
$$

 The last expression is finite iff the function $ z(\cdot) $ belongs to the space $ L_2(B,\nu_{\alpha,n}) $,
for instance, if it  is bounded. \par

 Moreover, it is interest to note that the variance of the r.v.  $  K_{n,B}(\alpha) \cdot z(\zeta),  $ i.e. the expression

$$
\sigma^2_{n,B}[z] := \Var \left\{ K_{n,B}(\alpha) \cdot z(\zeta) \right\} = K^2_{n,B}(\alpha) \cdot  ||z||^2L_2(\nu_{\alpha,n}) -
J^2_{\alpha,n}[z]
$$
  tends very  rapidly to zero as $  n \to \infty,  $  if for example $  \sup_k \alpha_k < 1.  $ \par
 Thus, we can use as in the second section by computation of this integral and by error estimation  the classical Monte -  Carlo
 method with application of Central Limit Theorem (CLT). \par
  There are not difficulties also to generate the multivariate ball beta  distribution as well as the polygonal beta distribution
generating.\par

\vspace{4mm}

\section{Parametric Volterra's integrals with weak singularities}

\vspace{4mm}

 Let $ z(\cdot) $ be  numerical function which dependent not only on the variable $ s, \ s \in S(n)  $ but on some
variable $  \theta, \ \theta \in \Theta, $  where $  \Theta  $ is any compact metrizable topological space:
$  z = z(s,\theta). $ We consider the following multidimensional parametric integral

$$
Q(\theta) =   \int_{S(n)} z(s,\theta) \  R_{\alpha,S}(s) \ ds =
$$

$$
 K_{n,S}(\alpha) \cdot \int_{S(n)} z(s,\theta) h_{\alpha}(s) \ ds = K_{n,S}(\alpha) \cdot {\bf E} z(\kappa, \theta), \eqno(5.1)
$$
where as before the r.v. $ \kappa $ has the polygonal beta distribution with parameters $  \alpha, n. $ \par
We offer for the parametric integral (5.1) computation the so-called "depending trial method", see \cite{Frolov1},
\cite{Grigorjeva1}, \cite{Ostrovsky1}, chapter 5, section 11:

$$
Q_N(\theta) := N^{-1} \sum_{i=1}^N z(\kappa_i, \theta), \eqno(5.2)
$$
where $ \{ \kappa_i \}  $ are independent polygonal distributed random vectors.\par
 In order to estimate a random uniform norm error for approximation \\ $ \sup_{\theta} |Q_N(\theta) - Q(\theta)|, $ we need to use the
Central Limit Theorem CLT on the Banach space of continuous functions $  C(\Theta) $  with ordinary uniform norm

$$
||f|| = \sup_{\theta \in \Theta} |f(\theta)|,
$$
 see \cite{Jain1},  \cite{Kozachenko1},  \cite{Ostrovsky1}, \cite{Dudley1},
\cite{Ostrovsky204}, \cite{Ostrovsky207}, \cite{Vaart1} etc. In detail, assume that the CLT in the space $ C(\Theta) $
for the sequence of the random fields $ \{ Q_N(\theta) - Q(\theta)   \}  $ there holds; then

$$
\lim_{N \to \infty} {\bf P} \left( \sqrt{N} \sup_{\theta \in \Theta} |Q_N(\theta) - Q(\theta)|  > u \right) =
{\bf P} \left( \sup_{\theta \in \Theta} |\xi(\theta)| > u   \right), \ u > 0, \eqno(5.3)
$$
where $ \xi(\theta) $ is mean zero continuous Gaussian random process (field) with at the same covariation function as the
random field $ z(\kappa_1,\theta) - {\bf E} z(\kappa_1,\theta):  $

$$
{\bf E} \xi(\theta_1) \xi(\theta_2) = K^2_{n,S}(\alpha) \int_{S(n)} z(s,\theta_1) z(s,\theta_2) \mu_{\alpha,n}(ds) - Q(\theta_1)Q(\theta_2).
\eqno(5.4)
$$

 The equality (5.3) allow us to construct the confidence interval in the uniform norm for calculated  function $  Q(\theta), $ see e.g.
 \cite{Ostrovsky1}, chapter 5, section 11; we investigate in the rest of this report the CLT in the considered case. \par

\vspace{3mm}

 Some notations.

$$
Y(s) := \sup_{\theta \in \Theta} |z(s,\theta)|, \ \rho(\theta_1,\theta_2) := \sup_{s \in S(n)} \left[ \frac{|z(s,\theta_1) - z(s,\theta_2)|}{Y(s)} \right].
$$

 The function $ (\theta_1, \theta_2) \to \rho(\theta_1,\theta_2) $ is bounded: $ \rho(\theta_1,\theta_2) \le 2 $
 continuous pseudo - metric on the set $  \Theta. $ We denote as usually by $ H(\Theta, \rho, \epsilon)  $ an entropy of the set $  \Theta $
 relative the semi-distance $  \rho $ at the point $ \epsilon,  $ i.e. the natural logarithm the minimal numbers of $ \rho \ -  $ balls
 of radii $  \epsilon, \ \epsilon > 0 $ which cover all the set $  \Theta. $ \par

\vspace{3mm}

{\bf Theorem.} \par
{\bf A.} {\it  Suppose  }

$$
{\bf E} |Y(\kappa)| =
\int_{S(n)} \sup_{\theta \in \Theta} |z(s,\theta)| \ \mu_{\alpha,n}(ds) < \infty. \eqno(5.5)
$$
{\it  Then the sequence  $ Q_N(\theta)  $ converges to $  Q(\theta) $ as $ N \to \infty  $ uniformly with probability one: }

$$
{\bf P} \left( \lim_{N \to \infty} \sup_{\theta \in \Theta} |Q_N(\theta) - Q(\theta)| \to 0   \right) = 1. \eqno(5.6)
$$

\vspace{3mm}

{\bf B.} {\it If }

$$
\int_{S(n)}  \sup_{\theta \in \Theta} |z(s,\theta) |^2 \ \mu_{\alpha,n}(ds) < \infty \eqno(5.7)
$$
{\it and}

$$
\int_0^1 H^{1/2}(\Theta, \rho, \epsilon) \ d \epsilon  < \infty, \eqno(5.8)
$$
{\it then the sequence of r.f. $ \sqrt{N} \left[ Q_N(\theta) - Q(\theta) \right]  $ satisfies the CLT in the space} $ C(\Theta,\rho). $\par

\vspace{3mm}

{\bf Remark 5.1.} The last condition (5.8) is satisfied if for example $  \Theta $ is closure of bounded open set in the space
$ R^d $ and $ \rho(\theta_1, \theta_2) \le C \ |\theta_1 - \theta_2|^{\gamma}, \ \gamma = \const \in (0,1];  $ in this case

$$
H(\Theta, \rho, \epsilon) \le C_2 + \  (d/\gamma) \ |\ln \epsilon |, \ \epsilon \in (0,1).
$$

 Moreover, this condition is satisfied if $ \dim_{\rho} \Theta < \infty.  $\par

\vspace{3mm}

{\bf Proof.} Let us  consider the  centered random field

$$
\lambda(\theta) = z(\kappa,\theta) - Q(\theta).
$$
This field  belongs to the (separable) Banach space $ C(\Theta,\rho)  $ with probability one and

$$
{\bf E} ||\lambda|| \le 2  \int_{S(n)} \sup_{\theta \in \Theta} |z(s,\theta)| \ \mu_{\alpha,n}(ds) < \infty.
$$
The first proposition of theorem follows from the well-known LLN in Banach spaces, theorem of Forte-Mourier. \par

 In order to  prove the second assertion {\bf B }, we need to use the famous result belonging  to Jain and Marcus \cite{Jain1}
about CLT in the space of continuous functions.  We have based of definition of the distance $ \rho $

$$
|\lambda(\theta_1) - \lambda(\theta_2)| \le 2 Y(\kappa)\rho(\theta_1,\theta_2).
$$
 Since $ {\bf E} Y^2(\kappa) < \infty  $ and  $ \int_0^1 H^{1/2}(\Theta, \rho, \epsilon) \ d \epsilon  < \infty, $
we deduce  needed for us CLT in the space $ C(\Theta, \rho).  $  \par

\vspace{3mm}

 Analogously may be considered the multiple parametric weak singular integral relative the ball beta distribution.\par

\end{document}